\documentclass{ifacconf}

\usepackage{natbib}

\usepackage{booktabs}

\usepackage{graphicx} 
\usepackage{amsmath} 
\usepackage{amssymb}  

\newcommand{\T}{^\mathrm{T}}
\newcommand{\Real}{\mathbb{R}}

\newcommand{\bio}{\sf{ }} 

\newtheorem{proposition}{Proposition}

\begin{document}

\begin{frontmatter}
\title{Optimality in cellular storage via the Pontryagin Maximum Principle}
\author[KUL]{Steffen Waldherr} and
\author[OVGU]{Henning Lindhorst}

\address[KUL]{KU Leuven, Department of Chemical Engineering,
Leuven, Belgium \\
email: steffen.waldherr@kuleuven.be}

\address[OVGU]{Institute for Automation Engineering,
Otto-von-Guericke-Universit\"at Magdeburg,
Magdeburg, Germany}

\begin{abstract}
We study an optimal control problem arising from a resource allocation problem in cellular metabolism.
A minimalistic model that describes the production of enzymatic vs.\ non-enzymatic biomass components from a single nutrient source is introduced.
The basic growth modes with this model are linear growth, where only the non-enzymatic component is produced, and exponential growth with only enzymatic components being produced.
Using Pontryagin's maximum principle, we derive the optimal growth trajectory depending on the model's parameters.
It turns out that depending on the parameters, either a single growth mode is optimal, or otherwise the optimal solution is a concatenation of exponential growth with linear growth.
Importantly, on the short time scale, the choice of growth mode depends only on catalytic rate constants and biomass weights of the two component types, whereas on longer time scales, where the nutrient amount becomes limiting, also the yield coefficients play a role.
\end{abstract}
\end{frontmatter}

\section{Introduction}

The regulation of cellular metabolism is the key step for microorganisms to ensure growth under dynamic environmental nutrient availability.
Due to the mechanisms of evolution, the regulation is adapted in order to optimize the cell's evolutionary fitness.
This optimality principle has given rise to a wide-spread use of optimization-based models to describe cellular metabolism.
The common approach is to formulate biophysical constraints on metabolic fluxes, such as thermodynamic or enzyme capacity constraints, combine them with a biologically relevant objective function, and solve the resulting optimization problem for the variables describing the cellular metabolism, such as metabolic reaction fluxes.
Experimental measurements confirm that actual metabolic states and dynamics often come close to the predictions of appropriately formulated optimization-based models \citep{VarmaPal1994,MahadevanEdw2002,SchuetzZam2012}.

A particular perspective on optimization-based models of cellular metabolism is to view them as a resource allocation problem.
In that case, the optimization problem is formulated to predict how cells would allocate building blocks and energy obtained from nutrients to different parts of the cellular machinery.
Previous studies from this perspective focussed on the allocation of resources to the nutrient uptake machinery or the proliferative machinery \citep{vandenBerg:1998bg,Giordano:2014wc} or the allocation of resources for different nutrient uptake variants \citep{vandenBerg:2002wo}.
These studies formulated simple, abstract models of the allocation problem, where the optimization problem could be solved analytically as an optimal control problem.
In addition, resource allocation problems for cellular metabolism have been studied with medium to large scale network models.
A steady state problem for optimal enzyme expression levels was proposed by \cite{Goelzer:2011vs} and solved for a \textit{Bacillus subtilis} model.
This model formalism was later extended by \cite{WaldherrOya2014} to dynamic conditions, the \emph{deFBA} method.
In applications of the deFBA method, it was shown that the resulting optimization problem provides a rationale for common observations in microorganisms, such as diauxic switches between different nutrients, overflow metabolism and reutilization of excreted substances, and adaptation to nutrient depletion \citep{WaldherrOya2014}.
Experimental measurements confirmed that enzyme expression levels qualitatively match optimization results from such models under a wide range of growth conditions \citep{Brien:2013fl}.

These previous resource allocation studies focused on the enzymatic machinery of the cell, because this machinery is required to drive metabolism and cell growth.
However, microorganisms also accumulate storage compounds.
In various bacteria, nitrogen limitation for example leads to the production of the polymer poly-$\beta$-hydroxybutyrate (PHB) as a carbon reserve, which is relevant for the biotechnological production of plastics \citep{Boe:1990jx,Beun:2002ea,FranzSon2011}.
Such storage compounds are difficult to account for in classic optimization-based models, as they are often related to a change in biomass composition and require to take dynamic processes over a certain time horizon into account.

Our numerical experiments with the deFBA method have resulted in the interesting observation that under some conditions, on short time horizons or with a low nutrient amount, a simple biomass objective function may actually yield the accumulation of storage compounds in the optimization result.
Intuitively, this seems to stem from the property that, compared to proteins as enzymatic compounds, storage compounds are usually faster to produce and have a higher yield due to lower specific energy requirements for producing them.
Based on this observation, our goal for this paper is to analytically characterize the conditions under which accumulation of storage compounds is optimal for maximizing a biomass objective.
Similar to previous theoretical resource allocation studies, we rely on an abstract model, that only distinguishes between allocation of resources to either storage compounds or enzymatic compounds.
Both compounds contribute to biomass, but only enzymatic compounds increase the metabolic capacity of the cells.
This model is simple enough to derive analytical results via the Pontryagin maximum principle.

\section{Model}
\label{sec:model}

We consider a simplistic metabolic-genetic network which consists of a nutrient metabolite $\bio{N}$, an energy metabolite $\bio{A}$, a storage macromolecule $\bio{M}$, an an enzyme $\bio{E}$.
The reaction stoichiometry is shown in Table~\ref{tab:network}.
\begin{table}
  \centering
  \caption{Resource allocation network}
  \label{tab:network}
  \begin{tabular}{lllr}
    \toprule
    Reaction & Stoichiometry & Enzyme & $k_\mathrm{cat}$ \\ \midrule
    $v_A$ & ${\bio{N}} \rightarrow {\bio{A}}$ & ${\bio{E}}$ & $\tilde k_A$ \\
    $v_M$ & $b_M\, {\bio{N}} + \tilde a_M b_M\, {\bio{A}} \rightarrow {\bio{M}}$ & ${\bio{E}}$ & $\tilde k_M$ \\
    $v_E$ & $b_E\, {\bio{N}} + \tilde a_E b_E\, {\bio{A}} \rightarrow {\bio{E}}$ & ${\bio{E}}$ & $\tilde k_E$ \\
    \bottomrule
  \end{tabular}
\end{table}

The network is reduced to three states by putting the quasi steady state assumption $\dot x_A = 0$ on the energy metabolite $\bio{A}$ \citep{WaldherrOya2014}.
From the relevant molar balance, this yields the algebraic constraint
\begin{equation}
  \label{eq:energy-qssa}
  v_A = \tilde a_M b_M v_M + \tilde a_E b_E v_E.
\end{equation}
Using the state vector $x = (x_N, x_M, x_E)\T$ for nutrient amount, storage molecule amount, and enzyme amount, respectively, yields the dynamics
\begin{equation}
  \label{eq:molar-balance-dynamics}
  \begin{aligned}
    \dot x_N &= - v_A - b_M v_M - b_E v_E = - a_M b_M v_M - a_E b_E v_E  \\
    \dot x_M &= v_M \\
    \dot x_E &= v_E,
  \end{aligned}
\end{equation}
with $a_M = \tilde a_M + 1$ and $a_E = \tilde a_E + 1$.
Since all fluxes are catalyzed by the single enzyme $\bio{E}$, the rates can be written as
\begin{equation}
  \label{eq:rate-transformation}
  v_i = u_i x_E,
\end{equation}
where $u_i$ is proportional to the amount of enzyme allocated to reaction $i \in \{A, M, E\}$.

Based on the quasi steady state reduction~\eqref{eq:energy-qssa}, we consider the enzyme allocation vector
\begin{equation}
  \label{eq:control-vector}
  u\T = (u_M, u_E)
\end{equation}
and the stoichiometric matrix
\begin{equation}
  \label{eq:stoich-matrix}
  S =
  \begin{pmatrix}
    -a_M b_M & - a_E b_E \\
    1 & 0 \\
    0 & 1
  \end{pmatrix},
\end{equation}
which together give rise to the network dynamics $\dot x = S u x_E$.

The key constraint on the control is the enzyme capacity constraint as considered in \cite{Goelzer:2011vs} and \cite{WaldherrOya2014}.
For the network considered here, this constraint is given by
\begin{equation}
  \label{eq:enzyme-capacity-constraint-full}
  \frac{v_A}{\tilde k_A} + \frac{v_M}{\tilde k_M} + \frac{v_E}{\tilde k_E} \leq x_E,
\end{equation}
where $\tilde k_i$ is the enzyme's catalytic constant for the reaction $v_i$.
From the quasi steady state reduction~\eqref{eq:energy-qssa} and the reaction rate transformation~\eqref{eq:rate-transformation},
we obtain a reduced constraint as
\begin{equation}
  \label{eq:enzyme-capacity-constraint}
  \frac{u_M}{k_M} + \frac{u_E}{k_E} \leq 1,
\end{equation}
with
\begin{equation*}
  \begin{aligned}
    \frac{1}{k_M} &= \frac{\tilde a_M b_M}{\tilde k_M} + \frac{1}{\tilde k_M} \\
    \frac{1}{k_E} &= \frac{\tilde a_E b_E}{\tilde k_E} + \frac{1}{\tilde k_E}.
  \end{aligned}
\end{equation*}

The cellular objective is typically related to biomass accumulation.
For the network considered here, the biomass is proportional to the amount of nutrients that is put into storage molecule or enzymes.
Any nutrient that is consumed for energy production does not contribute to biomass accumulation.
Using the biomass vector $b\T = (0, b_M, b_E)$, the biomass is thus given by
\begin{equation}
  \label{eq:biomass}
  b\T x = b_M x_M + b_E x_E.
\end{equation}

From the deFBA model as introduced in \cite{WaldherrOya2014}, the dynamic optimization problem
\begin{equation}
  \label{eq:optimization-problem}
  \begin{aligned}
    \max_{u(\cdot),x(\cdot)} &\ \int\limits_{0}^{T} b\T x(t) \;dt \\
    \textnormal{s.t. } & \dot{x} = S u x_E & \qquad \frac{u_M}{k_M} + \frac{u_E}{k_E} &\leq 1 \\
    & x(0) = x_0 & \qquad x,\ u &\geq 0 \\
  \end{aligned}
\end{equation}
is obtained.

\section{Background in optimal control}
\label{sec:background}

The analysis in this paper relies on necessary conditions for solutions of an optimal control problem of the general form
\begin{equation}
  \label{eq:oc-problem}
  \begin{aligned}
    \max_{u(\cdot),x(\cdot)} &\ \int\limits_{0}^{T} K(x,u) \;dt \\
    \textnormal{s.t. } & \dot{x} = f(x,u)  & \qquad g(u) &\geq 0 \\
    & x(0) = x_0 & \qquad h(x) &\geq 0.
  \end{aligned}
\end{equation}
We are using necessary conditions from Pontryagin's maximum principle with state constraints as summarized by \cite{Hartl:1995bv}.
Following their notation, all multipliers (greek letters) are defined as row vectors.

The Hamiltonian function is defined as
\begin{equation}
  \label{eq:hamiltonian}
  H(x,u,\lambda_0,\lambda) = \lambda_0 K(x,u) + \lambda f(x,u),
\end{equation}
and the Lagrangian function as
\begin{equation}
  \label{eq:lagrangian}
  L(x,u,\lambda_0,\lambda,\mu,\nu) = H(x,u,\lambda_0,\lambda) + \mu g(u) + \nu h(x).
\end{equation}
To simplify notation, functions which are evaluated on an optimal trajectory are written as
\begin{equation}
  \label{eq:2}
  H^\ast[t] =  H(x^\ast(t),u^\ast(t),\lambda_0,\lambda(t)),
\end{equation}
and equivalently for any other functions.

Let $\Omega = \{ u | g(u) \geq 0 \}$.
 Let $u^\ast(\cdot)$ and $x^\ast(\cdot)$ be an optimal solution of the control problem~\eqref{eq:oc-problem}, where $x^\ast(\cdot)$ has only finitely many junction times.
Then there exist a constant $\lambda_0 \geq 0$, a costate trajectory $\lambda(\cdot)$, and piecewise continuous multiplier functions $\mu(\cdot) \geq 0$ and $\nu(\cdot) \geq 0$ that satisfy:

\begin{enumerate}
\item The pointwise optimality condition
  \begin{equation}
    \label{eq:pointwise-optimality}
    u^\ast(t) = \arg \max_{u \in \Omega} H(x^\ast(t),u,\lambda_0,\lambda(t)) \\
  \end{equation}
\item The costate dynamics
  \begin{equation}
    \label{eq:costate-dynamics}
    \dot\lambda(t) = - \frac{\partial L^\ast}{\partial x}[t]
  \end{equation}
\item The first order and multiplier optimality conditions 
\begin{equation}
  \label{eq:multiplier-optimality-conditions}
  \begin{aligned}
  \frac{\partial L^\ast}{\partial u}[t] = 0 \\
  \mu(t) g^\ast[t] = 0 \\
  \nu(t) h^\ast[t] = 0.
  \end{aligned}
\end{equation}
\end{enumerate}

At the terminal time $t_f$, the condition
\begin{equation}
  \label{eq:terminal-condition}
  \lambda(t_f^-) = \gamma \frac{\partial h^\ast}{\partial x}[t_f],
\end{equation}
holds with 
\begin{equation}
  \label{eq:terminal-condition-2}
  \gamma \geq 0 \qquad \textnormal{and} \qquad \gamma h^\ast[t_f] = 0.
\end{equation}

The multipliers are required to satisfy
\begin{equation}
  \label{eq:multiplier-nonzero}
  (\lambda_0,\lambda(t),\mu(t),\nu(t),\gamma) \neq 0
\end{equation}
for all times $t$, i.e., at least one element, usually $\lambda_0$, must be non-zero.


\section{Results}
\label{sec:results}

\subsection{Optimization problem and candidate optimal solutions}
\label{sec:cand-optim-solut}

The optimization problem formulated in~\eqref{eq:optimization-problem} is non-linear.
However, by a suitable transformation of variables as in \cite{Jabarivelisdeh:2016wa}, it is equivalent to a linear optimization problem, for which the conditions of the maximum principle are necessary and sufficient for optimality \citep{Bressan:2007un}.
Therefore, we can use the Pontryagin maximum principle as necessary and sufficient conditions also for the problem as formulated in~\eqref{eq:optimization-problem}.

Thereby, the control constraints are
\begin{equation}
  \label{eq:control-constraints}
  g(u) =
  \begin{pmatrix}
    u_M \\
    u_E \\
    1 - \frac{u_M}{k_M} - \frac{u_E}{k_E},
  \end{pmatrix}
  \geq 0
\end{equation}
and the state constraints are
\begin{equation}
  \label{eq:state-constraints}
  h(x) = x \geq 0.
\end{equation}

The Hamiltonian~\eqref{eq:hamiltonian} is obtained as
\begin{equation}
  \label{eq:hamiltonian}
  H(x,u,\lambda_0,\lambda) = \lambda_0 b\T x + \lambda S u x_E
\end{equation}
and the Lagrangian~\eqref{eq:lagrangian} as
\begin{equation}
  \label{eq:lagrangian}
  L(x,u,\lambda_0,\lambda,\mu,\nu) = \lambda_0 b\T x + \lambda S u x_E + \mu g(u) + \nu h(x).
\end{equation}

Based on previous numerical analyses \citep{WaldherrOya2014,Lindhorst:2016vt}, we hypothesize that optimal solutions will be composed by up to three different phases:
\begin{itemize}
\item An \emph{exponential} phase where only enzyme is being produced:
  \begin{equation}
    \label{eq:control-exponential}
    u^\ast(t) = (0, k_E)\T.
  \end{equation}
  Within a time interval $t \in [ \tau_0, \tau_1]$, this yields the solution
  \begin{equation}
    \label{eq:state-exponential}
    x^\ast(t) =
    \begin{pmatrix}
      x_N(\tau_0) - a_E b_E x_E(\tau_0) \bigl( e^{k_E(t - \tau_0)} - 1 \bigr) \\
      x_M(\tau_0) \\
      x_E(\tau_0) e^{k_E(t-\tau_0)}
    \end{pmatrix}.
  \end{equation}
  Due to the state constraint $x_N(t) \geq 0$, the exponential solution can be maintained for a maximum time range of
  \begin{equation}
    \label{eq:exponential-max-time}
    \tau_1 - \tau_0 \leq \frac{1}{k_E} \ln \bigr( 1 + \frac{x_N(\tau_0)}{a_E b_E x_E(\tau_0)} \bigl).
  \end{equation}

\item A \emph{linear} phase where only the storage molecule is being produced:
  \begin{equation}
    \label{eq:control-linear}
    u^\ast(t) = (k_M, 0)\T.
  \end{equation}
  Within a time interval $t \in [ \tau_0, \tau_1]$, this yields the solution
  \begin{equation}
    \label{eq:state-linear}
    x^\ast(t) =
    \begin{pmatrix}
      x_N(\tau_0) - a_M b_M k_M x_E(\tau_0) (t - \tau_0) \\
      x_M(\tau_0) + k_M x_E(\tau_0) (t - \tau_0) \\
      x_E(\tau_0)
    \end{pmatrix}.
  \end{equation}
  Due to the state constraint $x_N(t) \geq 0$, the linear solution can be maintained for a maximum time range of
  \begin{equation}
    \label{eq:linear-max-time}
    \tau_1 - \tau_0 \leq \frac{x_N(\tau_0)}{a_M b_M k_M x_E(\tau_0)}.
  \end{equation}

\item A \emph{stationary} phase with zero flux, $u^\ast(t) = 0$, which yields the solution $x^\ast(t) = x^\ast(\tau_0)$ and can be maintained indefinitely.
\end{itemize}

\subsection{Costate dynamics}
\label{sec:costate-dynamics}

From~\eqref{eq:costate-dynamics}, the dynamics for the costate $\lambda$ are given by the differential equation
\begin{equation}
  \label{eq:costate-dynamics-general}
  \begin{aligned}
    \dot \lambda &= - \lambda_0 b\T - \lambda (0,\ 0,\ S u^\ast) - \nu
  \end{aligned}
\end{equation}

For each of the three phases defined above, we compute the dynamics and their solution for the costate $\lambda(t)$ within an interval $t \in [\tau_0, \tau_1]$, based on a terminal condition $\lambda(\tau_1)$.

\subsubsection{Exponential phase}
\label{sec:exponential-phase}

During exponential growth, we have $x_N(t) > 0$ and $x_E(t) > 0$.
From~\eqref{eq:multiplier-optimality-conditions}, this implies
\begin{equation}
  \nu_1 = 0 \qquad \nu_3 = 0.
\end{equation}
If $x_M(\tau_0) > 0$, then $\nu_2 = 0$, else $\nu_2$ can take any non-negative value.

With $u^\ast$ from~\eqref{eq:control-exponential} we have the costate dynamics
\begin{align}
  \label{eq:costate-dynamics-exponential}
    \dot \lambda_1 &= 0 \\
    \dot \lambda_2 &= - \lambda_0 b_M - \nu_2 \\
    \dot \lambda_3 &= - \lambda_0 b_E + a_E b_E k_E \lambda_1 - k_E \lambda_3,
\end{align}
yielding the solution
\begin{align}
  \label{eq:costate-exponential-lam1}
  \lambda_1(t) &= \lambda_1(\tau_1) \\
  \label{eq:costate-exponential-lam2}
  \lambda_2(t) &= \lambda_2(\tau_1) + \lambda_0 b_M (\tau_1 - t) + \int_t^{\tau_1} \nu_2(s) ds \\
  \label{eq:costate-exponential-lam3}
  \lambda_3(t) &= e^{k_E(\tau_1 - t)} \bigl(\lambda_3(\tau_1) - a_E b_E \lambda_1(\tau_1) + \lambda_0 \frac{b_E}{k_E}\bigr) \nonumber \\
               & \quad + a_E b_E \lambda_1(\tau_1) - \lambda_0 \frac{b_E}{k_E}.
\end{align}

\subsubsection{Linear phase}
\label{sec:linear-phase}

During linear growth, we have $x_i(t) > 0$ for all $i \in \{N,M,E\}$.
From~\eqref{eq:multiplier-optimality-conditions}, this implies that the multiplier for the state constraints is $\nu = 0$.

With $u^\ast$ from~\eqref{eq:control-linear} we have the costate dynamics
\begin{align}
  \label{eq:costate-dynamics-linear}
    \dot \lambda_1 &= 0 \\
    \dot \lambda_2 &= - \lambda_0 b_M \\
    \dot \lambda_3 &= - \lambda_0 b_E + a_M b_M k_M \lambda_1 - k_M \lambda_2,
\end{align}
yielding the solution
\begin{align}
  \label{eq:costate-linear-1}
  \lambda_1(t) &= \lambda_1(\tau_1) \\
  \label{eq:costate-linear-2}
  \lambda_2(t) &= \lambda_2(\tau_1) + \lambda_0 b_M (\tau_1 - t) \\
  \label{eq:costate-linear-3}
  \lambda_3(t) &= \lambda_3(\tau_1) + \frac{1}{2} \lambda_0 b_M k_M (\tau_1 - t)^2 \nonumber \\
               & \quad + (\lambda_0 b_E - a_M b_M k_M \lambda_1(\tau_1) + k_M \lambda_2(\tau_1)) (\tau_1 - t).
\end{align}

\subsubsection{Stationary phase}
\label{sec:stationary-phase}

With $u^\ast = 0$ we have the costate dynamics
\begin{equation}
  \label{eq:costate-dynamics-stationary}
  \dot \lambda = - \lambda_0 (b\T + \nu),
\end{equation}
yielding the solution
\begin{equation}
  \label{eq:costate-stationary}
  \lambda(t) = \lambda(\tau_1) + \lambda_0 b\T (\tau_1 - t) + \lambda_0 \int_t^{\tau_1} \nu(s) ds
\end{equation}
with
\begin{equation}
  \nu(t) x^\ast(t) = 0.
\end{equation}

\subsection{Optimality conditions}
\label{sec:optim-cond}

The analysis of optimality for the candidate solutions is based on the pointwise optimality condition~\eqref{eq:pointwise-optimality}.
For the system analyzed here, the pointwise optimality condition becomes:
\begin{equation}
  \label{eq:pointwise-optimality-condition}
  \begin{aligned}
  u^\ast(t) &= \arg \max_{u \in \Omega} x_E^\ast(t) \lambda S u \\
  &= \arg \max_{u \in \Omega} (\lambda_2 - a_M b_M \lambda_1) u_M + (\lambda_3 - a_E b_E \lambda_1) u_E
  \end{aligned}
\end{equation}
with $\Omega = \{ u \in \Real^2 : g(u) \geq 0 \}$, where we drop the argument $t$ for $\lambda$.

Singular arcs could exist if the coefficients of $u_M$ or $u_E$ in \eqref{eq:pointwise-optimality-condition} were constantly equal to zero along the arc, 
but this is inconsistent with the costate dynamics.
Singular arcs where the coefficients are equal and positive, but not necessarily constant, can not be ruled out at this point.

We consider optimal solutions that are composed of arcs representing the linear, exponential, or stationary solution.
For these arcs, the pointwise optimality conditions are specified as follows:
\begin{itemize}
\item Linear growth is pointwise optimal, if and only if the coefficient of $u_M$ in~\eqref{eq:pointwise-optimality-condition} is positive and larger than the coefficient of $u_E$, i.e.,
  \begin{equation}
    \label{eq:lin-condition}
    \lambda_2(t) > a_M b_M \lambda_1(t)
  \end{equation}
  and
  \begin{equation}
    \label{eq:lin-pointwise}
    -a_M b_M k_M \lambda_1(t) + k_M \lambda_2(t) > -a_E b_E k_E \lambda_1(t) + k_E \lambda_3(t).
  \end{equation}
\item Exponential growth is pointwise optimal, if and only if
  \begin{equation}
    \label{eq:exp-condition}
    \lambda_3(t) > a_E b_E \lambda_1(t)
  \end{equation}
  and
  \begin{equation}
    \label{eq:exp-pointwise}
   -a_E b_E k_E \lambda_1(t) + k_E \lambda_3(t) > -a_M b_M k_M \lambda_1(t) + k_M \lambda_2(t).
  \end{equation}
\item Stationary growth is optimal, if and only if
  \begin{equation}
    \label{eq:stationary-pointwise-1}
    \lambda_2(t) < a_M b_M \lambda_1(t)
  \end{equation}
  and
  \begin{equation}
    \label{eq:stationary-pointwise-2}
    \lambda_3(t) < a_E b_E \lambda_1(t).
  \end{equation}
\end{itemize}


\subsection{Optimal solutions on short horizons}
\label{sec:optim-solut-short}

In this section, we study optimal solutions to \eqref{eq:optimization-problem} with short terminal times $T$, i.e., such that the nutrient $\bio{N}$ will not be depleted within that time.
We prove that, depending on the model parameters, one of three growth dynamics is optimal on such short horizons:
Exponential growth~\eqref{eq:control-exponential}, linear growth~\eqref{eq:control-linear}, or exponential growth followed by a switch to linear growth.

From~\eqref{eq:terminal-condition} and~\eqref{eq:terminal-condition-2}, the terminal costate for short horizons is determined as $\lambda(T) = 0$.

\begin{proposition}[Exponential growth]
  The exponential control~\eqref{eq:control-exponential} is optimal for \eqref{eq:optimization-problem}, if and only if
  \begin{equation}
    \label{eq:condition-short-exponential}
    k_E b_E > k_M b_M
  \end{equation}
  and
  \begin{equation}
    \label{eq:time-exponential}
    x_N(0) > a_E b_E x_E(0) (e^{k_E T} - 1).
  \end{equation}
\end{proposition}

\begin{pf}
  We show pointwise optimality of an exponential solution arc with $\tau_0 = 0$ and $\tau_1 = T$.
  Such an arc is only feasible if \eqref{eq:time-exponential} holds.
  Independent of $x_M(0)$, the multiplier $\nu_2(s)$ in \eqref{eq:costate-exponential-lam2} can be set to $0$, because $\lambda_2(t)$ should be as small as possible for exponential growth to be optimal.
The first condition \eqref{eq:exp-condition} for pointwise optimality becomes
\begin{equation*}
  e^{k_E (T - t)} b_E - b_E > 0,
\end{equation*}
which is satisfied for $t < T$.
The second condition~\eqref{eq:exp-pointwise} is equivalent to
  \begin{equation*}
    \frac{e^{k_E(T-t)} - 1}{T - t} > \frac{b_M k_M}{b_E},
  \end{equation*}
which needs to hold for all $t \in (0,T)$.
Observing that $(e^{kx} -1)/x > k$ for $x > 0$ and $\lim_{x\rightarrow 0} (e^{kx}-1)/x = k$, this is equivalent to~\eqref{eq:condition-short-exponential}.
\qed 
\end{pf}

\begin{proposition}[Linear growth]
\label{prop:linear}
  The linear control~\eqref{eq:control-linear} is optimal for \eqref{eq:optimization-problem}, if and only if
  \begin{equation}
    \label{eq:time-short-linear}
    T < \frac{2(k_M b_M - k_E b_E)}{b_M k_M k_E}
  \end{equation}
and
\begin{equation}
  \label{eq:nutrient-short-linear}
  x_N(0) > T a_M b_M k_M x_E(0).
\end{equation}
\end{proposition}

Note that $T$ being positive means that   
\begin{equation}
  k_E b_E < k_M b_M
\end{equation}
needs to hold for linear growth to be optimal.

\begin{pf}
  We show pointwise optimality of a linear solution arc with $\tau_0 = 0$ and $\tau_1 = T$.
  Such an arc is only feasible if \eqref{eq:nutrient-short-linear} holds.
  The first condition~\eqref{eq:lin-condition} for pointwise optimality becomes
  \begin{equation*}
    b_M (T - t) > 0,
  \end{equation*}
  which is clearly satisfied for all $t$ between $0$ and $T$.
The second condition~\eqref{eq:lin-pointwise} becomes
\begin{equation*}
  k_M b_M (T - t) > \frac{1}{2} b_M k_M k_E (T - t)^2 + b_E k_E (T - t)
\end{equation*}
for all $t$ between $0$ and $T$, which is equivalent to~\eqref{eq:time-short-linear}.
\qed 
\end{pf}

\begin{proposition}[Exponential-linear growth]
  The optimal control is a switched function given by
  \begin{equation}
    \label{eq:exponential-linear-control}
    u^\ast(t) = \left\{
      \begin{aligned}
        \begin{pmatrix}
          0 \\ k_E
        \end{pmatrix} &\quad\textnormal{ for } t < \tau_1 \\
        \begin{pmatrix}
          k_M \\ 0
        \end{pmatrix} &\quad\textnormal{ for } \tau_1 < t < T
      \end{aligned}
      \right.
  \end{equation}
with
\begin{equation}
  \label{eq:switching-time-exponential-linear}
  \tau_1 = T - \frac{2(k_M b_M - k_E b_E)}{b_M k_M k_E},
\end{equation}
if and only if
\begin{equation}
  \label{eq:condition-short-explin}
  k_E b_E < k_M b_M,
\end{equation}
  \begin{equation}
    \label{eq:time-short-explin}
    T > \frac{2(k_M b_M - k_E b_E)}{b_M k_M k_E},
  \end{equation}
  and
  \begin{equation}
    \label{eq:nutrient-short-explin}
    x_N(0) > x_E(0) (a_E b_E (e^{k_E \tau_1} - 1) + a_M b_M e^{k_E \tau_1} (T - \tau_1)).
  \end{equation}
\end{proposition}

\begin{pf}
The control~\eqref{eq:exponential-linear-control} is feasible if~\eqref{eq:condition-short-explin}--\eqref{eq:nutrient-short-explin} are satisfied.
We then show pointwise optimality first for the linear arc active from $\tau_1$ to $T$, and then for the exponential arc active from $0$ to $\tau_1$.

For the linear arc, condition~\eqref{eq:lin-condition} is satisfied for $t < T$.
From~\eqref{eq:lin-pointwise}, we get
\begin{equation*}
  k_M b_M (T - t) > \frac{1}{2} b_M k_M k_E (T - t)^2 + b_E k_E (T - t)
\end{equation*}
for $t$ between $\tau_1$ and $T$
as in the proof of Proposition~\ref{prop:linear}, which is satisfied with~\eqref{eq:condition-short-explin} and~\eqref{eq:time-short-explin}.

For the exponential arc, we first compute the costate $\lambda(\tau_1)$ at the end of the exponential arc, which is the same as the beginning of the linear arc, from~\eqref{eq:costate-linear-1}--\eqref{eq:costate-linear-3}, obtaining
\begin{equation*}
  \begin{aligned}
    \lambda_1(\tau_1) &= 0 \\
    \lambda_2(\tau_1) &= \lambda_0 \frac{2(k_M b_M - k_E b_E)}{k_M k_E} \\
    \lambda_3(\tau_1) &= \frac{1}{2} \lambda_0 b_M k_M (T - \tau_1)^2 + \lambda_0 b_E (T - \tau_1).
  \end{aligned}
\end{equation*}
The first pointwise optimality condition~\eqref{eq:exp-condition} is satisfied because $\lambda_3 > 0$ and $\lambda_1 = 0$ for $t$ between $0$ and $\tau_1$.
The second pointwise optimality condition~\eqref{eq:exp-pointwise} becomes
\begin{equation}
  \label{eq:exponential-linear-inequality}
  \begin{aligned}
    e^{k_E (\tau_1 - t)} (\frac{1}{2} b_M k_M k_E (T - \tau_1)^2 + b_E k_E (T - \tau_1) + b_E) > \\
    b_M k_M (T - \tau_1) + b_M k_M (\tau_1 - t) + b_E.
  \end{aligned}
\end{equation}
Using~\eqref{eq:switching-time-exponential-linear} we get the relation
\begin{equation*}
  \frac{1}{2} b_M k_M k_E (T - \tau_1)^2 + b_E k_E (T - \tau_1) = b_M k_M (T - \tau_1),
\end{equation*}
which can be used to rewrite inequality~\eqref{eq:exponential-linear-inequality} to
\begin{equation*}
    \frac{e^{k_E (\tau_1 - t)} - 1}{\tau_1 - t} (b_M k_M (T - \tau_1) + b_E) > b_M k_M
\end{equation*}
for all $t$ between $0$ and $\tau_1$, which is equivalent to~\eqref{eq:condition-short-explin}.
\qed 
\end{pf}

\subsection{Optimal solutions on long horizons}
\label{sec:optim-solut-long}

In this section, we study optimal solutions to \eqref{eq:optimization-problem} with longer terminal times $T$, i.e., such that the nutrient $x_N$ will be depleted within the considered time range.
We consider three candidate optimal solutions for this problem.
These are composed by the exponential, linear, and exponential-linear solutions from the previous section, each followed by a stationary phase when the nutrient has been depleted.

Due to the depletion of the nutrient, $x_N(T) = 0$, the terminal costate $\lambda(T)$ will typically be non-zero.
We introduce a parameter $\gamma_1$, to be determined later, such that
\begin{equation}
  \label{eq:terminal-costate-stationary}
  \lambda_1(T) = \gamma_1 \geq 0,
\end{equation}
while $\lambda_2(T) = \lambda_3(T) = 0$.

Since the stationary arc is the same for all solution candidates, we start with analyzing optimality conditions for this.
Let $\tau_s$ be the time at which the stationary arc starts.
Based on the pointwise optimality conditions~\eqref{eq:stationary-pointwise-1} and~\eqref{eq:stationary-pointwise-2}, we obtain
\begin{equation}
  \label{eq:long-stationary-pointwise-1}
  \gamma_1 + \int_{\tau_s}^T \nu_1(s) ds \geq \lambda_0 \max \Bigl\{\frac{T - \tau_s}{a_M}, \frac{T - \tau_s}{a_E} \Bigr\}.
\end{equation}

\begin{proposition}[Linear-stationary growth]~ \\
\label{prop:lin-stat-growth}
  The optimal control is a switched function given by
  \begin{equation}
    \label{eq:linear-stationary-control}
    u^\ast(t) = \left\{
      \begin{aligned}
        \begin{pmatrix}
          k_M \\ 0
        \end{pmatrix} &\quad\textnormal{ for } t < \tau_s \\
        0 &\quad\textnormal{ for } t > \tau_s
      \end{aligned}
      \right.
  \end{equation}
with
\begin{equation}
  \label{eq:switching-time-linear-stationary}
  \tau_s = \frac{x_N(0)}{a_M b_M k_M x_E(0)}
\end{equation}
if and only if
\begin{equation}
  \label{eq:condition-linstat}
  a_E \geq a_M,
\end{equation}
  \begin{equation}
    \label{eq:time-linstat}
    (\frac{a_E}{a_M} - 1) b_E k_E T \geq \frac{1}{2}b_M k_M k_E \tau_s^2 + (\frac{a_E}{a_M} b_E k_E - b_M k_M) \tau_s,
  \end{equation}
  and
  \begin{equation}
    \label{eq:nutrient-linstat}
    x_N(0) < T a_M b_M k_M x_E(0).
  \end{equation}
\end{proposition}

\begin{pf}
Condition~\eqref{eq:nutrient-linstat} ensures that $x_N(t) = 0$ during the stationary phase, which is necessary for optimality.

From the condition~\eqref{eq:lin-condition} on piecewise optimality of the linear phase, we can derive
\begin{equation*}
\lambda_1(\tau_s) + \int_{\tau_s}^T \nu_1(s) ds \leq \frac{T - \tau_s}{\lambda_0 a_M},
\end{equation*}
which, together with the optimality condition~\eqref{eq:long-stationary-pointwise-1} on the stationary phase, implies~\eqref{eq:condition-linstat}.
We can then choose
  \begin{equation*}
    \gamma_1 = \lambda_0 \frac{T - \tau_s}{a_M},
  \end{equation*}
and $\nu_1 = 0$.

The costate at the end of the linear phase is taken from the start of the stationary phase.
Together with the costate dynamics for the linear phase in~\eqref{eq:costate-linear-1}--\eqref{eq:costate-linear-3}, this gives
\begin{equation}
\label{eq:costate-long-linear}
  \begin{aligned}
    \lambda_1(t) &= \gamma_1 \\
    \lambda_2(t) &= \lambda_0 b_M (T - t) \\
    \lambda_3(t) &= \lambda_0 \bigl( b_E (T - \tau_s) + \frac{1}{2} b_M k_M (\tau_s - t)^2 + ( b_E  \\
    & \quad \mbox{} - a_M b_M k_M \frac{\gamma_1}{\lambda_0} + b_M k_M (T - \tau_s)) (\tau_s - t) \bigr)
  \end{aligned}
\end{equation}
during the linear phase.
The first optimality condition~\eqref{eq:lin-condition} is satisfied for $t$ between $0$ and $\tau_s$ because of~\eqref{eq:condition-linstat}.
The second optimality condition~\eqref{eq:lin-pointwise} is rewritten as
\begin{equation}
  \label{eq:full-time-condition-linstat}
  \begin{aligned}
  -\frac{1}{2} b_M k_M k_E (\tau_s - t)^2 + (b_M k_M - b_E k_E) (\tau_s - t) \\
 + (\frac{a_E}{a_M} - 1) b_E k_E (T - \tau_s) > 0,\\
  \end{aligned}
\end{equation}
which needs to hold for all $t$ between $0$ and $\tau_s$.
Because the condition is convex in $(\tau_s - t)$, we only need to verify it at the endpoints $t = \tau_s$ and $t = 0$.
At $t = \tau_s$, it is satisfied with~\eqref{eq:condition-linstat} and~\eqref{eq:nutrient-linstat}.
Evaluating~\eqref{eq:full-time-condition-linstat} at $t = 0$ yields~\eqref{eq:time-linstat}.
\qed 
\end{pf}

\begin{proposition}[Exponential-stationary growth]
\label{prop:exp-stat}
If \begin{equation}
  \label{eq:condition-expstat}
  a_M \geq a_E,
\end{equation}
\begin{equation}
  \label{eq:max-time-expstat}
  x_N(0) < (e^{k_E T} - 1) a_E b_E x_E(0),
\end{equation}
and either
\begin{equation}
  \label{eq:expstat-condition-1}
  b_E k_E \geq b_M k_M,
\end{equation}
or, alternatively to~\eqref{eq:expstat-condition-1},
\begin{equation}
  \label{eq:expstat-condition-2}
  \tau_s \leq (1 - \frac{a_E}{a_M}) T,
\end{equation}
where 
\begin{equation}
  \label{eq:switching-time-expstat}
  \tau_s = \frac{1}{k_E} \ln(1 + \frac{x_N(0)}{a_E b_E x_E(0)}),
\end{equation}
then the optimal control is given by the switched function
  \begin{equation}
    \label{eq:exponential-stationary-control}
    u^\ast(t) = \left\{
      \begin{aligned}
        \begin{pmatrix}
          0 \\ k_E
        \end{pmatrix} &\quad\textnormal{ for } t < \tau_s \\
        0 &\quad\textnormal{ for } t > \tau_s.
      \end{aligned}
      \right.
  \end{equation}
\end{proposition}

  \begin{pf}
Condition~\eqref{eq:max-time-expstat} ensures that $x_N(t) = 0$ during the stationary phase, which is necessary for optimality.

From the condition~\eqref{eq:exp-condition} on piecewise optimality of the exponential phase, we can derive
\begin{equation*}
\lambda_1(\tau_s) + \int_{\tau_s}^T \nu_1(s) ds \leq \frac{T - \tau_s}{a_M},
\end{equation*}
which, together with the optimality condition~\eqref{eq:long-stationary-pointwise-1} on the stationary phase, implies~\eqref{eq:condition-expstat}.
We can then choose
  \begin{equation*}
    \gamma_1 = \lambda_0 \frac{T - \tau_s}{a_E},
  \end{equation*}
and $\nu_1 = 0$.

From~\eqref{eq:costate-exponential-lam1}--\eqref{eq:costate-exponential-lam3}, the costate during the exponential phase is then given by
\begin{equation*}
  \begin{aligned}
  \lambda_1(t) &= \gamma_1 \\
  \lambda_2(t) &= \lambda_0 b_M (T - t) \\ 
  \lambda_3(t) &= e^{k_E(\tau_s - t)} (\lambda_0 b_E (T - \tau_s) - a_E b_E \gamma_1 + \lambda_0 \frac{b_E}{k_E}) \\
               & \quad + a_E b_E \gamma_1 - \lambda_0 \frac{b_E}{k_E}.
  \end{aligned}
\end{equation*}
The first optimality condition~\eqref{eq:exp-condition} is satisfied for $t$ between $0$ and $\tau_s$.
The second optimality condition~\eqref{eq:exp-pointwise} is rewritten as
\begin{equation}
  \label{eq:exp-pointwise-expstat}
  \begin{aligned}
    (e^{k_E (\tau_s - t)} - 1) \frac{b_E}{b_M k_M} > - (\frac{a_M}{a_E} - 1) (T - \tau_s) + (\tau_s - t)
  \end{aligned}
\end{equation}
which needs to hold for all $t$ between $0$ and $\tau_s$.
This is a condition about the intersection of an exponential with a linear function.
Without the interval constraint on $t$, it could be solved via Lambert's W-function, but for $t$ bounded by $\tau_s$, a closed form solution of this seems not to be achievable.

Conditions~\eqref{eq:expstat-condition-1} and~\eqref{eq:expstat-condition-2} in the proposition are derived from relaxations of~\eqref{eq:exp-pointwise-expstat} considering cases where
(a) the slope of the exponential on the left hand side is always larger than the slope of the linear function on the right hand side, or
(b) the right hand side is always negative, respectively.
\qed 
  \end{pf}

Note that Proposition~\ref{prop:exp-stat} gives only a sufficient condition on optimality of the exponential-stationary solution, which is not expected to be tight due to the approximation of~\eqref{eq:exp-pointwise-expstat} done in the proof.
Still, for specific parameter values, a necessary and sufficient condition can be obtained numerically by computing values for $t$ where~\eqref{eq:exp-pointwise-expstat} holds with equality using Lambert's W-function, and checking whether such a $t$ lies within $0$ and $\tau_s$.

\begin{proposition}[Exponential-linear-stationary growth]~ \\
\label{prop:exp-lin-stat-growth}
  The optimal control is a switched function given by
  \begin{equation}
    \label{eq:exp-linear-stationary-control}
    u^\ast(t) = \left\{
      \begin{aligned}
        \begin{pmatrix}
          0 \\ k_E
        \end{pmatrix} &\quad\textnormal{ for } t < \tau_1 \\
        \begin{pmatrix}
          k_M \\ 0
        \end{pmatrix} &\quad\textnormal{ for } \tau_1 < t < \tau_s \\
        0 &\quad\textnormal{ for } t > \tau_s,
      \end{aligned}
      \right.
  \end{equation}
if and only if
\begin{equation}
  \label{eq:condition-explinstat-1}
  a_E \geq a_M,
\end{equation}
and there exist $\tau_1$, $\tau_s$ with $0 < \tau_1 < \tau_s < T$ such that
  \begin{align}
  \label{eq:exp-lin-stat-time-1}
    x_N(0) - a_E b_E x_E(0) (e^{k_E \tau_1} - 1) - a_M b_M k_M (\tau_s - \tau_1) &= 0 \\
  \label{eq:exp-lin-stat-time-2}
    \frac{1}{2} b_M k_M k_E (\tau_s - \tau_1)^2 + (b_E k_E - b_M k_M) (\tau_s - \tau_1) &  \nonumber \\  + (1 - \frac{a_E}{a_M}) b_E k_E (T - \tau_s) &= 0,
  \end{align}
and
\begin{equation}
  \label{eq:condition-explinstat-2}
  b_E k_E \geq b_M k_M (1 - k_E (\tau_s - \tau_1)).
\end{equation}
\end{proposition}

\begin{pf}
We first consider optimality of the linear phase.
Condition~\eqref{eq:condition-explinstat-1} is derived as in the proof of Proposition~\ref{prop:lin-stat-growth}, also with the choice of
\begin{equation*}
  \gamma_1 = \lambda_0 \frac{T - \tau_s}{a_M}.
\end{equation*}

For the further steps, it is helpful to denote
\begin{equation*}
  \begin{aligned}
  z(t) &= \frac{1}{2} b_M k_M k_E (\tau_s - t)^2 + b_E k_E (\tau_s - t) \\
  &\qquad\mbox{} - (\frac{a_E}{a_M} - 1) b_E k_E (T - \tau_s).
  \end{aligned}
\end{equation*}
The costate for the linear phase is given as in~\eqref{eq:costate-long-linear}.
Condition~\eqref{eq:lin-pointwise} can then be rewritten as
\begin{equation}
\label{eq:convex-condition-explinstat}
  b_M k_M (\tau_s - t) - z(t) > 0
\end{equation}
which has to hold for $\tau_1 < t < \tau_s$.
Because the condition is convex, it needs only be evaluated at the end points of the interval.
At $t = \tau_s$, \eqref{eq:convex-condition-explinstat} is satisfied with~\eqref{eq:condition-explinstat-1} and $T > \tau_s$.

At $t = \tau_1$, \eqref{eq:convex-condition-explinstat} then becomes
\begin{equation}
  \label{eq:z-condition}
  z(\tau_1) \leq b_M k_M (\tau_s - \tau_1).
\end{equation}

We next consider the exponential phase.
The costate is obtained from~\eqref{eq:costate-linear-1}--\eqref{eq:costate-linear-3} and with the values from the linear phase as terminal condition at $\tau_1$, yielding
\begin{equation*}
  \begin{aligned}
  \lambda_1(t) &= \gamma_1 \\
  \lambda_2(t) &= \lambda_0 b_M (T - t) \\ 
  \lambda_3(t) &= \lambda_0 \bigl(e^{k_E(\tau_1 - t)} (\frac{z(\tau_1)}{k_E} + \frac{b_E}{k_E}) + \mbox{} \\
  &\qquad \frac{a_E}{a_M} b_E (T - \tau_s) - \frac{b_E}{k_E} \bigl).
  \end{aligned}
\end{equation*}
Condition~\eqref{eq:exp-condition} is then rewritten as
\begin{equation*}
  e^{k_E(\tau_1 - t)}(z(\tau_1) + b_E) > b_E,
\end{equation*}
which is satisfied for $t \in (0,\tau_1)$ with $\tau_1 > 0$ and $z \geq 0$.
Condition~\eqref{eq:exp-pointwise} becomes
\begin{equation}
  \label{eq:exp-lin-stat-exp-condition}
  e^{k_E(\tau_1 - t)}(z(\tau_1) + b_E) > b_M k_M (\tau_s - T) + b_E,
\end{equation}
which, together with~\eqref{eq:z-condition} implies that we have to choose $z(\tau_1) = b_M k_M (\tau_s - T)$, which is stated in~\eqref{eq:exp-lin-stat-time-2}.
That means that~\eqref{eq:exp-lin-stat-exp-condition} becomes an equality at $t = \tau_1$.
In order for the inequality~\eqref{eq:exp-lin-stat-exp-condition} to be satisfied for $t < \tau_1$, the slope of the exponential on the left side has to be steeper than the slope of the linear function on the right side, leading to~\eqref{eq:condition-explinstat-2} in Proposition~\ref{prop:exp-lin-stat-growth}.
\qed 
\end{pf}

\section{Conclusions}
\label{sec:conclusions}

The paper provides a comprehensive analysis of optimal solutions in a resource allocation problem, with a focus on the distinction between linear growth related to cellular storage mechanisms and the typical exponential growth obtained by investing into enzymatically active biomass.

The short time scale is relevant for optimization based models of cellular metabolism, which have previously been studied in \cite{WaldherrOya2014} and \cite{Lindhorst:2016vt}.
The key criterion for optimality is which growth mode is instantaneously faster for biomass accumulation, represented by the value of $b_E k_E$ for exponential growth and $b_M k_M$ for linear growth.
For actual biological systems, the production of enzymes is often slower and less efficient in mass accumulation than storage molecules.
In such a case, the exponential solution requires a certain time to ``overtake'' the linear solution, and this time should be considered in optimization based models where the goal would be to obtain the more realistic exponential solution.

On longer time horizons, where the nutrients are going to be depleted, obviously the yield coefficients for the different growth modes play a role as well.
In fact, which growth mode is optimal just before the depletion of the nutrient seems to depend only on the yield coefficients.
As production of storage molecules typically requires less energy and is more efficient compared to enzyme production, optimization based models for realistic biological systems should be expected to predict a linear growth phase or accumulation of storage molecules under nutrient limitation.
The question whether that would be ``optimal'' in an evolutionary sense is interesting, but it is beyond the scope of this study.


\begin{thebibliography}{16}
\providecommand{\natexlab}[1]{#1}
\providecommand{\url}[1]{\texttt{#1}}
\providecommand{\urlprefix}{URL }
\expandafter\ifx\csname urlstyle\endcsname\relax
  \providecommand{\doi}[1]{doi:\discretionary{}{}{}#1}\else
  \providecommand{\doi}{doi:\discretionary{}{}{}\begingroup
  \urlstyle{rm}\Url}\fi

\bibitem[{Beun et~al.(2002)Beun, Dircks, Van~Loosdrecht, and
  Heijnen}]{Beun:2002ea}
Beun, J.J., Dircks, K., Van~Loosdrecht, M.C.M., and Heijnen, J.J. (2002).
\newblock {Poly-$\beta$-hydroxybutyrate metabolism in dynamically fed mixed
  microbial cultures}.
\newblock \emph{Water Research}, 36(5), 1167--1180.

\bibitem[{Boe and Lovrien(1990)}]{Boe:1990jx}
Boe, I.N. and Lovrien, R.E. (1990).
\newblock {Energy reserves and storage polymers in intact bacteria analyzed by
  metabolic calorimetry}.
\newblock \emph{Thermochimica Acta}, 172, 115--122.

\bibitem[{Bressan and Piccoli(2007)}]{Bressan:2007un}
Bressan, A. and Piccoli, B. (2007).
\newblock \emph{{Introduction to the Mathematical Theory of Control}}.
\newblock American Institute of Mathematical Sciences.

\bibitem[{Franz et~al.(2011)Franz, Song, Ramkrishna, and Kienle}]{FranzSon2011}
Franz, A., Song, H.S., Ramkrishna, D., and Kienle, A. (2011).
\newblock {Experimental and theoretical analysis of
  poly($\beta$-hydroxybutyrate) formation and consumption in \emph{Ralstonia
  eutropha}}.
\newblock \emph{Biochem. Engin. J.}, 55(1), 49--58.

\bibitem[{Giordano et~al.(2014)Giordano, Mairet, Gouz{\'e}, Geiselmann, and
  de~Jong}]{Giordano:2014wc}
Giordano, N., Mairet, F., Gouz{\'e}, J.L., Geiselmann, J., and de~Jong, H.
  (2014).
\newblock {Dynamic optimisation of resource allocation in microorganisms}.
\newblock In \emph{21st International Symposium on Mathematical Theory of
  Networks and Systems}, 887--889.

\bibitem[{Goelzer et~al.(2011)Goelzer, Fromion, and Scorletti}]{Goelzer:2011vs}
Goelzer, A., Fromion, V., and Scorletti, G. (2011).
\newblock {Cell design in bacteria as a convex optimization problem}.
\newblock \emph{Automatica}, 47, 1210--1218.

\bibitem[{Hartl et~al.(1995)Hartl, Sethi, and Vickson}]{Hartl:1995bv}
Hartl, R.F., Sethi, S.P., and Vickson, R.G. (1995).
\newblock {A Survey of the Maximum Principles for Optimal Control Problems with
  State Constraints}.
\newblock \emph{SIAM Review}, 37(2), 181--218.

\bibitem[{Jabarivelisdeh and Waldherr(2016)}]{Jabarivelisdeh:2016wa}
Jabarivelisdeh, B. and Waldherr, S. (2016).
\newblock {Improving Bioprocess Productivity Using Constraint-Based Models in a
  Dynamic Optimization Scheme}.
\newblock In \emph{6th IFAC Symposium on Foundations of Systems Biology in
  Engineering}.

\bibitem[{Lindhorst et~al.(2016)Lindhorst, Lucia, Findeisen, and
  Waldherr}]{Lindhorst:2016vt}
Lindhorst, H., Lucia, S., Findeisen, R., and Waldherr, S. (2016).
\newblock {Modeling metabolic networks including gene expression and
  uncertainties}.
\newblock arXiv:1609.08961 [math.OC].

\bibitem[{Mahadevan et~al.(2002)Mahadevan, Edwards, and
  Doyle~III}]{MahadevanEdw2002}
Mahadevan, R., Edwards, J.S., and Doyle~III, F.J. (2002).
\newblock {Dynamic flux balance analysis of diauxic growth in Escherichia
  coli.}
\newblock \emph{Biophys.\ J.}, 83(3), 1331--1340.

\bibitem[{O'Brien et~al.(2013)O'Brien, Lerman, Chang, Hyduke, and
  Palsson}]{Brien:2013fl}
O'Brien, E.J., Lerman, J.A., Chang, R.L., Hyduke, D.R., and Palsson, B.O.
  (2013).
\newblock {Genome-scale models of metabolism and gene expression extend and
  refine growth phenotype prediction}.
\newblock \emph{Molecular Systems Biology}, 9(1), 693.

\bibitem[{Schuetz et~al.(2012)Schuetz, Zamboni, Zampieri, Heinemann, and
  Sauer}]{SchuetzZam2012}
Schuetz, R., Zamboni, N., Zampieri, M., Heinemann, M., and Sauer, U. (2012).
\newblock {Multidimensional optimality of microbial metabolism.}
\newblock \emph{Science}, 336(6081), 601--604.

\bibitem[{van~den Berg et~al.(1998)van~den Berg, Kiselev, Kooijman, and
  Orlov}]{vandenBerg:1998bg}
van~den Berg, H.A., Kiselev, Y.N., Kooijman, S.A.L.M., and Orlov, M.V. (1998).
\newblock {Optimal allocation between nutrient uptake and growth in a microbial
  trichome}.
\newblock \emph{Journal of Mathematical Biology}, 37(1), 28--48.

\bibitem[{van~den Berg et~al.(2002)van~den Berg, Kiselev, and
  Orlov}]{vandenBerg:2002wo}
van~den Berg, H.A., Kiselev, Y.N., and Orlov, M.V. (2002).
\newblock {Optimal allocation of building blocks between nutrient uptake
  systems in a microbe.}
\newblock \emph{Journal of Mathematical Biology}, 44(3), 276--296.

\bibitem[{Varma and Palsson(1994)}]{VarmaPal1994}
Varma, A. and Palsson, B.{\O}. (1994).
\newblock {Stoichiometric flux balance models quantitatively predict growth and
  metabolic by-product secretion in wild-type Escherichia coli W3110.}
\newblock \emph{Appl Environ Microbiol}, 60(10), 3724--3731.

\bibitem[{Waldherr et~al.(2015)Waldherr, Oyarz{\'u}n, and
  Bockmayr}]{WaldherrOya2014}
Waldherr, S., Oyarz{\'u}n, D.A., and Bockmayr, A. (2015).
\newblock {Dynamic optimization of metabolic networks coupled with gene
  expression}.
\newblock \emph{Journal of Theoretical Biology}, 365, 469--485.

\end{thebibliography}
\end{document}